\documentclass[10pt]{amsart}
\usepackage{amsthm,amssymb}
\usepackage{graphicx}

\newtheorem{lemma}{Lemma}
\newtheorem{theorem}{Theorem}
\newtheorem{definition}{\rm D e f i n i t i o n}

\theoremstyle{definition}

\newtheorem{example}{\rm E x a m p l e}

\newtheorem{remark}{\rm R e m a r k}

\newcommand{\supp}{\mbox{supp}}
\newcommand{\ordin}{\mbox{ord}}

\newcommand{\cl}{\mbox{cl}}

\newcommand{\loc}{\mbox{loc}}
\newcommand{\var}{\mbox{var}}

\newcommand{\dyn}{D}
\newcommand{\us}{U}

\begin{document}

\title[Distributions with dynamic test functions]{Distributions with dynamic test functions and 
multiplication by discontinuous functions}

\author[V.Derr] {V.~Derr}

\address{Faculty of Mathematics, Udmurtia State University \\ Universitetskya St., 1 (building 4), Izhevsk, 426034, Russia}

\email{derr@uni.udm.ru}

\author[D.Kinzebulatov] {D.~Kinzebulatov}

\address{Department of Mathematics and Statistics, University of Calgary \\ 2500 University Drive N.W., Calgary, Alberta, Canada T2N 1N4}

\email{knzbltv@udm.net}



\begin{abstract}
As follows from the Schwartz Impossibility Theorem, multiplication of two distributions is in general impossible.
Nevertheless, often one needs to multiply a distribution by a discontinuous function, not by an arbitrary distribution. In the present paper we construct a space of distributions where the general operation of multiplication by a discontinuous function is defined, continuous, commutative, associative and for which the Leibniz product rule holds. In the new space of distributions, the classical delta-function $\delta_\tau$ extends to a family of delta-functions $\delta_\tau^\alpha$, dependent on the \textit{shape} $\alpha$.
We show that the various known definitions of the product of the Heaviside function and the  delta-function in the classical space of distributions $\mathcal D'$ become particular cases of the multiplication in the new space of distributions, and provide the applications of the new space of distributions to the ordinary differential equations which arise in optimal control theory. 
Also, we
compare our approach of the Schwartz distribution theory with the approach of the
Colombeau generalized functions algebra, where the general operation of multiplication of two distributions is defined.

\end{abstract}

\maketitle

\section{Introduction}

The theory of distributions created by L.~Schwartz \cite{Scw} in the 1950s is an important part of modern mathematics.
However, in 1962 R.~Courant pointed out certain insufficiencies of the Schwartz distributions as generalizations of functions, in particular, the absence of a general multiplication of distributions \cite{Cour}. This multiplication is important for the definition of a solution for a wide class of ordinary and partial differential equations \cite{Col1}. It was mentioned  by J.F.Colombeau in \cite{Col1} that there is a strong belief by many mathematicians that a general multiplication of distributions is impossible in the Schwartz distribution theory, unless we drop essential properties (such as continuity), which is not acceptable for applications to differential equations.
%

Along with that, in many cases in the theory of differential equations one needs to be able to multiply a distribution not by an arbitrary distribution, but by a discontinuous function only (which is also a distribution). For instance, the product of a distribution and a discontinuous function arises in ordinary differential equations and partial differential equations of optimal control theory, e.g. see \cite{Ramp,Silv,Mil,Der3}, where a distribution represents an impulsive control. 
In our paper we construct a new space of distributions which are extensions of the classical distributions from the space of continuous test functions to the space of \textit{dynamic test functions}. In this space the
operation of multiplication of a distribution by a discontinuous function is defined, commutative, associative and continuous. 

Further we provide some applications of a new space of distributions to the theory of ordinary differential equations. 
Also, let us mention that since the operation of multiplication of a distribution and a discontinuous function can be induced by the operation of multiplication of two distributions, its properties allow us to compare, in a certain sense, various definitions of the product of two distributions in the classical space of distributions with continuous test functions $\mathcal D'$, which were proposed, in particular, in \cite{Bag1,Bag2,Sar,Sar2,Sar3}, in attempt to obtain a better definition.

The operation of multiplication of a distribution by a discontinuous function (as a particular case of a general operation of multiplication of two distributions) is defined in the Colombeau generalized functions algebra, which was applied successfully to many problems related to ordinary and partial differential equations \cite{Col1,ColM1,ColM2,ColM3, Col3, Col2,Col4}.
In the present paper we show that the commutative, associative and continuous operation of multiplication of a distribution by a discontinuous function can be also defined in the classical setting of the Schwartz distribution theory. Below we
compare our approach based on the Schwartz distribution theory with the approach of the
Colombeau generalized functions algebra (Section 4).




\subsection{Continuity properties of the operation of multiplication}

Let us describe some natural requirements imposed on continuity properties of the operation of multiplication.
Below the term ``regular distribution'' refers to a distribution which corresponds to a locally-summable function.

The operation of multiplication of a distribution by a discontinuous function has to have the following property: \\ [3mm]
\textbf{Property $\mathcal A$}.                                             
\textit{If $g$ is a discontinuous function, $v$ is a distribution, $\{v_n\}_{n=1}^\infty$ is a sequence of regular distributions, $v_n \to v$,
then $gv_n \to gv$. If $v$ is a regular distribution, then $gv$ 
is regular and corresponds to the usual product of $g$ and $v$}. \\ [-2mm]

A space of distributions where the operation of multiplication by a discontinuous function has property $\mathcal A$ was constructed in \cite{DerKin4}. 
This operation of multiplication arises from the following ordinary differential equation \cite{Mil}, \\ [-4mm]
\begin{equation}
\label{intro_eq1}
\dot{x}=f(t,x)+g(t,u(t))v
\end{equation}
where $v$ is a distribution (impulsive control), $f$, $g$ are continuously differentiable, $u$ is a control, which is in general discontinuous function \cite{Mil}. 

However, property $\mathcal A$ is too weak for consideration in the following ordinary differential equation with distributions which arises in many problems of optimal control theory \cite{Ses,Silv,Ramp}, \\ [-4mm]
\begin{equation}
\label{intro_eq2}
\dot{x}=f(t,x)+g(t,x)v, 
\end{equation}
where $v$ is a distribution, $f$, $g$ are continuously differentiable. 
In the general case, if the distribution $v$ is not regular, it is natural to expect that a solution $x$ is discontinuous \cite{Ses,Silv,Ramp,Mil}. Thus, 
the product of a distribution and a discontinuous function in (\ref{intro_eq2}) must have the following property.  \\ [1mm]
\textbf{Property $\mathcal B$}.                                                                                                                                                                                   
\textit{If $g$ is a discontinuous function, $v$ a distribution, $\{v_n\}_{n=1}^\infty$ a sequence of regular distributions, $v_n \to v$, and $\{g_n\}_{n=1}^\infty$ 
is a sequence of continuous functions, $g_n \to g$,
then $g_nv_n \to gv$. If $v$ is a regular distribution, then 
$gv$ is regular and corresponds to the usual product of $g$ and $v$}. \\ [-2mm]

Thus, in order to define the continuous operation of multiplication of distributions by discontinuous functions possessing property $\mathcal B$, which is important for the applications to ordinary differential equations, one needs to distinguish different approximations $\{v_n\}_{n=1}^\infty$ of the distribution $v$ in $\mathcal D'$, and different approximations $\{g_n\}_{n=1}^\infty$ of the discontinuous function $g$. 

This idea was formalized in \cite{ColM1,ColM2,ColM3} (see the references therein)
in the framework of the Colombeau generalized functions theory
(see more detailed discussion in Section 4) and, in less generality, in \cite{Mil2,Mil,Bre,Ses} (see further references therein),
where the consideration of the sequence $\{v_n\}_{n=1}^\infty$ as a part of the system (\ref{intro_eq2}) in order to provide the uniqueness of solution was proposed. 


In the present paper we formalize this idea in the classical setting of the Schwartz distribution theory, i.e., in the \textit{space of distributions with dynamic test functions}, where
the operation of multiplication by discontinuous function has property $\mathcal B$. 



\subsection{Distributions with dynamic test functions}
Let us briefly describe our approach. In distribution theory (as a theory of space $\mathcal D'$) the product of a distribution and a continuous function is defined by
\begin{equation}
\label{proddef}
(gf,\varphi)=(f,g\varphi)
\end{equation}
where $f$ is a distribution, $g$ is a continuous function and $\varphi$ is a continuous test function. If $g$ is discontinuous, then $g\varphi$ is no longer a test function, thus the value of the right-hand side of (\ref{proddef}) is undefined. In order to employ the definition (\ref{proddef}), we need to extend the space of test functions so that it contains discontinuous functions. In this case the value of the right-hand side of (\ref{proddef}) would be defined for discontinuous $g$ also. 

In the present paper we construct the space of distributions with the test functions possibly discontinuous.
We show that any classical distribution can be extended from the space of continuous test functions to our space of test functions. 
In particular, the classical delta-function $\delta_\tau$ extends into a family of delta-functions $\{\delta_\tau^\alpha\}$, where the parameter 
$\alpha$
such that 
$\int_{-1/2}^{1/2} \alpha(t)dt=1$, 
is called the \textit{shape} of a delta-function $\delta_\tau^\alpha$.
In the new space of distributions the product of the Heaviside function $\theta_\tau$ and the delta-function $\delta_\tau^\alpha$ is given by \\ [-2mm]
\begin{equation}
\label{proddef2}
\theta_\tau\delta_\tau^\alpha=\biggl(\int_0^{1/2}\alpha(t)dt \biggr)\delta_\tau^\beta, \\ [-1mm]
\end{equation}
where $\int_0^{1/2}\alpha(t)dt \in \bf R$ is a constant, and $\beta$ is another shape. Similarly, the family of delta-sequences extends into subfamilies of delta-sequences having different shapes. Then a result similar to (\ref{proddef2}) can be obtained if delta-functions in (\ref{proddef2}) are replaced by terms of delta-sequences having the same shapes.

We should mention that (\ref{proddef2}) generalizes known definitions of the product of the Heaviside function and the classical delta-function in $\mathcal D'$ given by
\begin{equation}
\label{stupidmult}
\theta_\tau\delta_\tau=c\delta_\tau,
\end{equation}
where $c \in \bf R$, which were proposed (in particular, through the multiplication of two distributions) in \cite{Kur,Fil,Tvr1,Ses,Sar,Sar2,Sar3}, in order to obtain the optimal (from a certain point of view) value of $c \in \bf R$, e.g. definition (\ref{stupidmult}) for $c=1/2$ satisfies the formal Leibniz rule for the differentiation of the product. 

In contrast to definitions of form (\ref{stupidmult}), the operation of multiplication in the new space of distributions is continuous, commutative and associative, the operation of multiplication in the new space of distributions satisfies the Leibniz product rule.


Apart of the problem of multiplication of distributions, our study is motivated by the problem of the correct definition of the solution of an ordinary differential equation with distributions. Further we provide an application of our space of distributions to nonlinear ordinary differential equations containing one delta-function $\delta_\tau^\alpha$, which arise in optimal control theory \cite{Ses,Ramp,Silv,Mil}. 

Note that the set of points of concentration of delta-functions in $v$ and even its cardinality, as well as the set of points of discontinuity of $x$, can not be fixed a priori, since in many optimal control problems it is a subject of optimization \cite{Silv,Ramp,Ses}. This leads to the necessity of considering the space of functions of bounded variation, and, as a result, to certain complications in contrast to the case when the set of points of discontinuity of a function, which is multiplied by a distribution, is finite and fixed. 

In order to construct this space of distributions
we introduce the supplementary notion of the \textit{dynamic function}. Most of the definitions for the ordinary functions (such as limit at point, support, boundedness, continuity, bounded variation) are transferred to dynamic functions without significant changes.
 
\section{Dynamic test functions \\ [3mm]} 

\subsection{Notations}
Let $I=(a,b) \subset \bf R$ be a fixed open interval, in general unbounded.
By $\textbf{R}^n$ and $\textbf{M}_n$ we denote the space of vectors
and the space of square matrices of order $n$ with real elements, respectively. We denote by $T(g) \subset I$
the set of points of discontinuity of a function $g:I \to \bf R$.

Consider the algebra of functions $g:I \to \bf R$ possessing one-sided limits 
\begin{equation*}
g(a+)=\lim_{t \to a+}g(t), \quad g(b-)=\lim_{t \to b-}g(t), \quad g(\tau+)=\lim_{t \to \tau+}g(t), \quad g(\tau-)=\lim_{t \to \tau-}g(t)
\end{equation*}
for any $\tau \in I$.
We identify functions having the same one-sided limits, and possibly having different values in their points of discontinuity.
We denote this algebra of functions (i.e., equivalence classes) by $\mathbb G=\mathbb G(I)$ and endow it with the norm
\begin{equation*}
\|g\|_{\mathbb G}=\sup_{t \in I}\max\{|g(t+)|,|g(t-)|\}.
\end{equation*}
With this norm $\mathbb G$ is a Banach algebra \cite{Der2}.
The algebra $\mathbb G$ is called the \textit{algebra of regulated functions} \cite{Dieu}.
The set of points of discontinuity $T(g)=\{t \in I: g(t+) \ne g(t-)\}$ of a regulated function $g \in \mathbb G$ is at most countable \cite{Dieu}.

For a function $g \in \mathbb{G}$ let us consider a partition $d=\{t_i\}_{i=1}^n \subset I$, 
$t_1<\dots<t_n$ of the interval $I$.
Denote 
\begin{equation*}
\var_{I}(g) \doteq \sup_d \sum_{i=1}^n|g(t_i-)-g(t_{i-1}+)| \geqslant 0
\end{equation*}
(\textit{total variation}). If $\var_{I}(g)<\infty$, we say that $g$ is a \textit{function of bounded variation} on $I$.
We denote this algebra by
$\mathbb{BV}=\mathbb{BV}(I)$, and define the norm
\begin{equation*}
\|g\|_{\mathbb{BV}}=|g(a+)|+\var_{I}(g) \\ [3mm].
\end{equation*}
It follows that $\mathbb{BV}$ is a Banach algebra.

We denote by $\mathbb C=\mathbb C(I)$ the subalgebra of continuous elements of $\mathbb G$. Denote by $\mathbb{CBV}$ the subalgebra of continuous elements of $\mathbb{BV}$.
Denote by $\mathbb L=\mathbb L(I)$ the Banach algebra of functions Lebesgue summable on $I$ with the integral norm.
Further $\mathbb{AC}=\mathbb{AC}(I)$ stands for the Banach algebra of absolutely continuous functions \cite{Dun}.
Let $\mathbb L(\loc)$ and $\mathbb{AC}(\loc)$ be the algebras of locally summable and locally absolutely continuous functions, respectively.   

We denote $\sigma_t(g) \doteq g(t+)-g(t-)$.

Let us denote by $\mathbb F=\mathbb F(I)$ the algebra of functions $I \to \mathbf R$ with the operations of point-wise addition, multiplication by an element of $\bf R$ and multiplication. The elements of $\mathbb F$ are called the \textit{ordinary} functions. 

The definitions given above are transferred without significant changes to the case of a finite closed interval.


Notice, that ignoring  values of discontinuous functions at their points of discontinuity leads to a further simplification of the space of distributions.

\subsection{Dynamic functions} In this subsection we introduce the supplementary notion of a \textit{dynamic function} which is needed for the construction of the space of test function.
It generalizes the notion of an ordinary function. Let $J=[-\frac{1}{2},\frac{1}{2}]$. 

\begin{definition} 
A map
$f:I \to \mathbb F(J)$
is called the dynamic function.
\end{definition}

We denote the value of a dynamic 
function $f$ at a point $t$ by $f(t)(\cdot)$. 

Denote the set of dynamic functions by $d\mathbb F=d\mathbb F(I)$.                     
We say that dynamic functions $f_1$, $f_2$ are equal, if $f_1(t)(\cdot)=f_2(t)(\cdot)$ in $\mathbb F(J)$ for all
$t \in I$. 
Also, the operations of addition, multiplication by an element of $\bf R$ and multiplication are defined in $d\mathbb F$ as follows:
for dynamic functions $f_1,f_2 \in d\mathbb F$ we define their sum $f_1+f_2 \in d\mathbb F$ by the formula
\begin{equation*}
(f_1+f_2)(t)(\cdot)=f_1(t)(\cdot)+f_2(t)(\cdot)
\end{equation*}
for all $t \in I$;
similarly, we define the product
\begin{equation*}
(f_1f_2)(t)(\cdot)=f_1(t)(\cdot)f_2(t)(\cdot);
\end{equation*}
for a given dynamic function $f \in d\mathbb F$ and $\lambda \in \bf R$ we define 
\begin{equation*}
(\lambda f)(t)(\cdot)=\lambda f(t)(\cdot) 
\end{equation*}
for all $t \in I$. Thus, $d\mathbb F$ form an algebra. 

We define the inclusion of the algebra of ordinary functions $\mathbb F$ into the algebra of dynamic functions $d\mathbb F$ by the map:
with an ordinary function $\hat{f} \in \mathbb F$ we associate a dynamic function $f \in d\mathbb F$ defined by
\begin{equation*}
f(t)(\cdot) \equiv \hat{f}(t)
\end{equation*}
for all $t \in I$.                             

For a given $t$, we call $f(t)(\cdot)$ a \textit{dynamic value} of $f$ at point $t$.
A dynamic value of $f$ which is identically equal to a constant function is called an \textit{ordinary value}. If $f \in d\mathbb F$ has an ordinary value at a point $t$,
we denote it by $f(t) \in \bf R$. 
Then a dynamic function is ordinary if and only if it has ordinary values at all points. 
We denote the set of points, where $f$ has ordinary values, by $\us(f) \subset I$. Let $\dyn(f) \doteq I \setminus \us(f)$. \\

Let $f \in d\mathbb F$, $f(t)=0$ ($t \ne \tau$), with dynamic value $f(\tau)(\cdot)=\gamma(\cdot)$. Construct a sequence of ordinary functions $\{f_n\}_{n=1}^\infty$ defined by
\begin{equation}
f_n(t)=\left\{
\begin{array}{l}
\gamma\bigl(n(t-\tau)\bigr), \text{ if } t \in \bigl(\tau-\frac{1}{2n},\tau+\frac{1}{2n}\bigr), \\ [2mm]
0, \quad \text{ otherwise} \\ [2mm]
\end{array}
\right.
\end{equation}
The sequence $\{f_n\}_{n=1}^\infty \subset \mathbb F$ is called the \textit{sequential representation} of $f \in d\mathbb{F}$.


We define the composition $g \circ f \in d\mathbb F$ of a dynamic function $f \in d\mathbb F$ and an ordinary function $g \in \mathbb F$ by
\begin{equation*}
(g \circ f)(t)(\cdot)=g \circ (f(t)(\cdot))
\end{equation*}
for all $t \in I$. 
Accordingly, define the absolute value of $f \in d\mathbb F$ by
\begin{equation*}
|f|(t)(\cdot) \doteq |f(t)(\cdot)|
\end{equation*}
for all $t \in I$. Also, we define the support \\
\begin{equation*}
\supp f=\cl\{t \in I: f(t)(\cdot) \text{ is not ordinary or } f(t) \ne 0\} \\ [2mm]
\end{equation*}
where $\cl$ stands for closure in $I$. For a given set $M \subset I$ let us define
\begin{equation*}
\sup_M f \doteq \sup_{t \in M} \sup_{s \in J} f(t)(s), \quad
\inf_M f \doteq \inf_{t \in I} \inf_{s \in J} f(t)(s). 
\end{equation*}
We say that $f$ is bounded on $M$, if $\sup_M|f|<\infty$.

A dynamic function $f$ is called nonnegative (positive) (denote $f \geqslant 0$ ($f>0$)) if $f(t)(s) \geqslant 0$ ($f(t)(s)>0$) for all $t \in I$.
A dynamic function $f$ is called nonpositive (negative) if $-f$ is nonnegative (positive).

Let $f \in d\mathbb F$ be bounded in a right neighborhood of a point $\tau \in I$.
We say that $f$ tends to $c \in \bf R$ from the right as $t \to \tau+$ if for any $\varepsilon>0$ there exists an $\eta>0$ such that
\begin{equation}
\sup_{(\tau,\tau+\eta)} |f-c|<\varepsilon.
\end{equation}
We denote the right-sided limit $c$ by $f(\tau+)$. The left-sided limit and the limit of a dynamic functions at a point are defined similarly. Analogously, define the one-sided limits at $a$, $b$.

We say that $f \in d\mathbb F$ is continuous at $\tau \in I$, if
\begin{equation}
f(\tau)=f(\tau+)=f(\tau-). 
\end{equation}

Notice that according to the notational convention above the notations used in this definition imply that $f$ has an ordinary value $f(\tau) \in \bf R$.
If $f$ is not continuous at  $\tau$, then we call $f$ discontinuous at this point.

Clearly, if in the definitions given above a dynamic function $f$ is an ordinary function, then these definitions coincide with the usual ones. 

\subsection{Special algebras of dynamic functions}
Motivated by our interest in ordinary differential equations, we consider the following subalgebras of the algebra of dynamic functions.

Denote by $d\mathbb G$ an algebra of dynamic functions $f$ having dynamic values $f(t)(\cdot) \in \mathbb G(J)$ and such that the
one-sided limits 
\begin{equation*}
f(a+), \quad f(b-), \quad f(\tau+), \quad f(\tau-)
\end{equation*}
exist for all $\tau \in I$. We call the elements of $d\mathbb G$ the \textit{dynamic regulated functions} (since $\mathbb G$ is a factor-space, $d\mathbb G$ is also a factor-space; definitions of the supremum and infimum are changed respectively; if there is no need, we will not make special remarks, and call the elements of $d\mathbb G$ the dynamic functions).
Let us define the norm in the algebra of dynamic regulated functions $d\mathbb G$ by
\begin{equation*}
\|f\|_{d\mathbb G}=\sup_{I}|f|.
\end{equation*}

We have the following properties of the elements of $d\mathbb G$.

\begin{lemma}
\label{lem1}
Any dynamic regulated function $f \in d\mathbb G$ has ordinary values on $I$ except for at most a countable set.
\end{lemma}

The proof of Lemma \ref{lem1} and other statements of Sections 2 and 3 are provided in Section 5.

Given $f \in d\mathbb G$, define a function $\hat{f}$ by $\hat{f}(t)=f(t)$ for any $t \in \us(f)$.
We call $\hat{f}$ the \textit{ordinary part} of $f$.
As follows from Lemma \ref{lem1}, the ordinary part $\hat{f}$ is defined everywhere on $I$ 
except for at most a countable set. Denote $\ordin(f) \doteq \hat{f}$.

\begin{lemma}
\label{lem2}
The ordinary part $\hat{f}$ of a dynamic regulated function $f \in d\mathbb G$ 
belongs to $\mathbb G$, and the equalities $f(t+)=\hat{f}(t+)$, $f(t-)=\hat{f}(t-)$ hold for all $t \in I$.
\end{lemma}


\begin{lemma}
\label{lem3}
The set of points of discontinuity $T(f)$ of any dynamic regulated function $f \in d\mathbb G$ is at most countable.
\end{lemma}

Let us define an inclusion of $\mathbb G$ into $d\mathbb G$. To do this we associate with any regulated function $\hat{f} \in \mathbb G$
a dynamic regulated function $f \in d\mathbb G$ having dynamic values
\begin{equation}
f(t)(s)=\hat{f}(\tau-) \text{ for } s \in [-1/2,0), \quad 
f(t)(s)=\hat{f}(\tau+) \text{ for } s \in (0,1/2] 
\end{equation}
for all $t \in I$. 

Let $f \in d\mathbb G$ be such that
\begin{equation*}
f(t)(\cdot) \in \mathbb{AC}(J), \quad f(t)(-1/2)=f(t-), \quad f(t)(1/2)=f(t+)
\end{equation*}
for all $t \in I$.
The set of such dynamic regulated functions with operations induced from $d\mathbb G$ form a
subalgebra $s\mathbb G \subset d\mathbb G$. Notice that $T(g)=\dyn(g)$ ($g \in s\mathbb G$).

We say that $f \in s\mathbb G$ is a \textit{dynamic function of bounded variation}, if $\ordin(f) \in \mathbb{BV}$ and
\begin{equation*}
\sum_{t \in \dyn(f)}\var_{s \in J}f(t)(s)<\infty, 
\end{equation*}
where $\dyn(f) \subset I$ is at most countable according to Lemma \ref{lem2}.
The set of dynamic functions of bounded variation with operations induced from $s\mathbb G$ forms a subalgebra.
The subalgebra of such dynamic functions is denoted by $s\mathbb{BV}$ and endowed with the norm
\begin{equation*}                          
\|f\|_{s\mathbb{BV}}=|f(a+)|+\|\hat{f}_c\|_{\mathbb{BV}}+\sum_{t \in \dyn(f)}\var_J (f(t)(\cdot)) \\ [1mm]
\end{equation*}
where $\hat{f}_c \in \mathbb{CBV}$ is a continuous part of  $\hat{f}=\ordin(f) \in \mathbb{BV}$ (the \textit{Jordan decomposition of functions of bounded variation}). 

\begin{example} The dynamic Heaviside function $\theta_\tau^\beta \in s\mathbb{BV}$ is defined by
\begin{equation*}
\theta_\tau^\beta(t)=\left\{
\begin{array}{l}
1, \quad t>\tau, \\
0, \quad t<\tau,
\end{array}
\right. \qquad
\theta_\tau^\beta(\tau)(\cdot)=\beta(\cdot), \\
\end{equation*}
where $\beta \in \mathbb{AC}(J)$ is such that $\beta(-1/2)=0$, $\beta(1/2)=1$.
\end{example}

\begin{example}
Since $\mathbb G \subset d\mathbb G$, the usual Heaviside function $\theta_\tau \in \mathbb G$ is in $d\mathbb G$. According to the definition of the inclusion of
$\mathbb G$ to $d\mathbb G$, we have that 
\begin{equation*}
\theta_\tau(t)=\left\{
\begin{array}{l}
0, \quad t<\tau, \\
1, \quad t>\tau,
\end{array}
\right. \qquad
\theta_\tau(\tau)(s)=\left\{
\begin{array}{l}
0, \quad s \in [-\frac{1}{2},0), \\
1, \quad s \in (0,\frac{1}{2}].
\end{array}
\right.
\end{equation*}
\end{example}


\subsection{Test functions}        
We denote by $\mathcal D$ the space of classical test functions, i.e., the space of the real-valued continuous functions having compact support in $I \subset \bf R$, which is endowed with the standard topology \cite{Shi}. 

We denote by $\mathcal T$ the space of elements $\varphi \in d\mathbb G$ having compact support $\supp\varphi \subset I$.
We say that $\{\varphi_n\}_{n=1}^\infty \subset \mathcal T$ converges to $\varphi \in \mathcal T$ if $\varphi_n \to \varphi$ in $d\mathbb G$ and 
there exists a closed interval $[c,d] \subset I$ such that $\supp\varphi_n \subset [c,d]$ for all $n=1,2,\dots$ Thus,
the classical space of continuous test functions $\mathcal D$ is contained in $\mathcal T$ as a subspace.

Notice that since for $\varphi \in \mathcal T$, $g \in d\mathbb G$ their product $g\varphi$ belongs to $\mathcal T$,
algebra $\mathcal T$ is an ideal in $d\mathbb G$. Further, if $\{\varphi_n\}_{n=1}^\infty \subset \mathcal T$ and
$\varphi_n \to \varphi$, then $g\varphi_n \to g\varphi$ in $\mathcal T$.

We call $\mathcal T$  the space of \textit{dynamic test functions}.
The proof of the following theorem is similar to the proof of an analogous statement for
$\mathcal D$ \cite{Ios}.  

\begin{theorem}
\label{topteo}
$\mathcal T$ is a locally convex topological vector space.
\end{theorem}
\section{Distributions with dynamic test functions \\ [3mm]}

Following standard notation \cite{Shi}, we denote by $\mathcal D'$ the space of distributions with continuous test functions, i.e., the space of continuous linear functionals
$\mathcal D \to \bf R$ \cite{Shi}.

We denote by $\mathcal T'$ the space of \textit{distributions with dynamic test functions}, i.e., the space of continuous linear functionals $\mathcal T \to \bf R$. 
The value of a distribution $f \in \mathcal T'$ on a test function $\varphi \in \mathcal T$ is denoted by $(f,\varphi) \in \bf R$.

\begin{example}
Let $f \in d\mathbb G$ (in particular, $f \in \mathbb G$). 
We define a distribution, which is also denoted by $f$, by the formula
\begin{equation}
\label{reg}
(f,\varphi)=\int_I\hat{f}(t)\hat{\varphi}(t)dt
\end{equation}        
where $\varphi \in \mathcal T$, $\hat{\varphi}=\ordin(\varphi)$,
$\hat{f}=\ordin(f)$.
The linearity and continuity of $f$ follows from the properties of the integral and the definition of 
the ordinary part of a dynamic regulated function; a distribution $f \in \mathcal T'$ is called a \textit{regular distribution}.
\end{example}

\begin{example}
Define a distribution $\delta_\tau^\alpha$ by the formula
\begin{equation}
\label{delta}
(\delta_\tau^\alpha,\varphi)=\int_J\varphi(\tau)(s)\alpha(s)ds
\end{equation}
where $\varphi \in \mathcal T$ and $\alpha\in \mathbb L(J)$ is such that
\begin{equation}
\label{normalization}
\int_J \alpha(s)ds=1.
\end{equation}                                 
The linearity and continuity of $\delta_\tau^\alpha$ follows from the properties of the integral and the definition of the convergence in $\mathcal T$, so
$\delta_\tau^\alpha \in \mathcal T'$. 
The distribution $\delta_\tau^\alpha \in \mathcal T'$ is called the \textit{delta-function} concentrated at the point $\tau \in I$ and having \textit{shape} $\alpha \in \mathbb L(J)$.
Now, for continuous test functions $\varphi \in \mathcal D$ 
\begin{equation*} (\delta_\tau^\alpha,\varphi)=\int_J\varphi(\tau)\alpha(s)ds=\varphi(\tau)\int_J\alpha(s)ds=\varphi(\tau).
\end{equation*}
Thus, the delta-function $\delta_\tau \in \mathcal D'$ extends to the family of delta-functions $\delta_\tau^\alpha \in \mathcal T'$. 
For a given delta-function $\delta_\tau^\alpha$ we construct a sequence $\{\omega_n^\alpha\}_{n=1}^\infty$,
\begin{equation*}
\omega^\alpha_n(t)=\left\{
\begin{array}{l}
n\alpha(n(t-\tau)), \text{ if } t \in (\tau-\frac{1}{2n},\tau+\frac{1}{2n}), \\ [3mm]
0, \text{ otherwise }
\end{array}
\right.
\end{equation*}
We call $\{\omega^\alpha_n\}_{n=1}^\infty$ the \textit{delta-sequence having shape} $\alpha \in \mathbb L(J)$.
\end{example}

\begin{example}
\label{ex_newdelta}
Let us define a distribution $\delta_\tau^\lambda \in \mathcal T'$, which is also called \textit{delta-function}, by 
\begin{equation}
\label{newdelta}
(\delta_\tau^\lambda,\varphi)=\lambda\varphi(\tau+)+(1-\lambda)\varphi(\tau-),
\end{equation}
where $\lambda \in \bf R$ and the restriction $\delta_\tau^\lambda|_{\mathcal D}=\delta_\tau \in \mathcal D'$.
Clearly, (\ref{newdelta}) can not be obtained from (\ref{delta}). The notion of a delta-sequence is not defined for the delta-function (\ref{newdelta}). 
\end{example}
            
\begin{remark}
\label{rm_newdelta}
In the general case, the delta-function in $\mathcal T'$ is defined as an affine combination of delta-functions (\ref{delta}) and (\ref{newdelta}). Below we consider mainly the delta-functions of the form (\ref{delta}).
\end{remark}              
                                                                                                       
In the space of distributions $\mathcal T'$ the linear operations of addition and multiplication by elements of $\bf R$ are
introduced in a standard way, so $\mathcal T'$ is a linear space. Let
$f_n \to f$ ($f_n,f \in \mathcal T'$) in $\mathcal T'$, if
$(f_n,\varphi) \to (f,\varphi)$ as $n \to \infty$ for any test function $\varphi \in \mathcal T$. The proof of the following lemma is similar to the proof of analogous theorem for $\mathcal D'$ \cite[p.65]{Shi}.

\begin{lemma} 
\label{lem4}                                         
If $\{f_n\}_{n=1}^\infty$ converges in $\mathcal T'$, and $\varphi_n \to 0$ in $\mathcal T$,
then $(f_n,\varphi) \to 0$.
\end{lemma}

\begin{theorem}
\label{limitteo}
Let the sequence $\{f_n\}_{n=1}^\infty \subset \mathcal T'$ be such that for any $\varphi \in \mathcal T$ the
sequence $\{(f_n,\varphi_n)\}_{n=1}^\infty$ converges as $n \to \infty$. The functional $f$ defined on $\mathcal T$ by 
\begin{equation}
\label{limit}
(f,\varphi)=\lim\limits_{n \to \infty} (f_n,\varphi)
\end{equation}
is linear and continuous and thus belongs to $\mathcal T'$.
\end{theorem}

\begin{theorem}
\label{ext_teo}
Any distribution in $\mathcal D'$ can be extended from $\mathcal D$ to $\mathcal T$.
\end{theorem}

Thus, the space $\mathcal T'$ is a space of extensions of distributions in $\mathcal D'$ from $\mathcal D$ to $\mathcal T$.

Notice, that the restriction of the distribution $\delta_\tau^\alpha-\delta_\tau^\gamma \in \mathcal T'$, where $\alpha$, $\gamma$ are the shapes of delta-functions, from $\mathcal T$ to $\mathcal D$, is a zero distribution in $\mathcal D'$. Hence, since the shapes $\alpha$, $\gamma$ can be chosen arbitrarily, 
addition of $\delta_\tau^\alpha-\delta_\tau^\gamma$ to an extension from $\mathcal D$ to $\mathcal T$ of any distribution in $\mathcal D'$, gives another extension.
Consequently, any distribution in $\mathcal D'$ has infinitely many extensions from $\mathcal D$ to $\mathcal T$.

A distribution $f \in \mathcal T'$ is called nonnegative (nonpositive) if for any $\varphi \in \mathcal T$ such that $\varphi \geq 0$ we have
$(f,\varphi) \geq 0$ (or $(f,\varphi) \leq 0$, respectively). 


Let us define in $\mathcal T'$ the operation of multiplication of distributions by the elements of $d\mathbb G$.
Since for any $g \in d\mathbb G$, $\varphi \in \mathcal T$ their product $g\varphi$ belongs to $\mathcal T$, we may define
\begin{equation}
\label{multdef}
(gf,\varphi) \doteq (f,g\varphi). \\ [2mm]
\end{equation}
where $f \in \mathcal T'$.
In particular, since $\mathbb G \subset d\mathbb G$,
any distribution in $\mathcal T'$ can be multiplied by any piece-wise continuous function.

\begin{theorem}
\label{cont_teo}
Let $f_n \to f$ in $\mathcal T'$, $g_n \to g$ in $d\mathbb G$. Then $g_nf_n \to gf$ in $\mathcal T'$.
\end{theorem}

The operation of multiplication in $\mathcal T'$ defined by (\ref{multdef}) is continuous, commutative and associative in the sense that
$(gh)f=g(hf)$ in $\mathcal T'$ for any $g,h \in d\mathbb G$, $f \in \mathcal T'$. 

\begin{example}
\label{ex21}
The product of the delta-function $\delta_\tau^\alpha$ and the Heaviside function $\theta_\tau^\beta \in s\mathbb G$ is given by the formula
\begin{equation}
\label{prod_w1}
(\theta_\tau^\beta\delta_\tau^\alpha,\varphi)=(\delta_\tau^\alpha,\theta_\tau^\beta\varphi)=\int_J\varphi(\tau)(s)\beta(s)\alpha(s)ds.
\end{equation}
If $\int_J\beta(s)\alpha(s)ds \ne 0$, then (\ref{prod_w1}) can be rewritten as
\begin{equation}
\label{prod33}
\theta_\tau^\beta\delta_\tau^\alpha=\biggl(\int_J\beta(s)\alpha(s)ds \biggr)\delta_\tau^\gamma,
\end{equation}
where 
$\gamma(\cdot)=\beta(\cdot)\alpha(\cdot)/\int_J\beta(s)\alpha(s)ds$ satisfies (\ref{normalization}). 
\end{example}

\begin{example}
\label{ex2}
The product of the delta-function $\delta_\tau^\alpha$ and the Heaviside function $\theta_\tau \in \mathbb G$ is given by the formula \\ [-5mm]
\begin{equation}
\label{prod_w2}
(\theta_\tau\delta_\tau^\alpha,\varphi)=(\delta_\tau^\alpha,\theta_\tau\varphi)=\int_{0}^{1/2}\varphi(\tau)(s)\alpha(s)ds.
\end{equation}
If $\int_{0}^{1/2}\alpha(r)dr \ne 0$, then (\ref{prod_w2}) can be rewritten as 
\begin{equation*}
\theta_\tau\delta_\tau^\alpha=\biggl(\int_0^{1/2}\alpha(s)ds\biggr)\delta_\tau^\gamma,
\end{equation*}
where the shape of the delta-function $\delta_\tau^\gamma$ is defined by $\gamma(s)=0$ for $s \in [-1/2,0)$, $\gamma(s)=\alpha(s)/\int_{0}^{1/2}\alpha(r)dr$ for $s \in (0,1/2]$.
\end{example}

\begin{remark}
\label{rem1}
Let us show that the equality (\ref{prod33}) can be obtained if $\delta_\tau^\alpha$ and $\theta_\tau^\beta$ are replaced by the terms of the corresponding delta-sequence $\{\omega^\alpha_n\}_{n=1}^\infty$ and sequential representation $\{f_n\}_{n=1}^\infty$, respectively. 
Let $\int_J\beta(s)\alpha(s)ds \ne 0$.
We have that
\begin{equation*}
f_n(t)\omega_n^\alpha(t)=\left\{
\begin{array}{l}
\beta(n(t-\tau))n\alpha(n(t-\tau)), \text{ if } t \in (\tau-\frac{1}{2n},\tau+\frac{1}{2n}), \\
0, \text{ otherwise},
\end{array}
\right.
\end{equation*}
i.e., \\ [-5mm]
\begin{equation*}
f_n(t)\omega_n^\alpha(t)=\int_J \beta(s)\alpha(s)ds\left\{
\begin{array}{l}
n\frac{\beta(n(t-\tau))\alpha(n(t-\tau))}{\int_J \beta(s)\alpha(s)ds}, \text{ if } t \in (\tau-\frac{1}{2n},\tau+\frac{1}{2n}), \\
0, \text{ otherwise},
\end{array}
\right.
\end{equation*}
where $\gamma(s)=\beta(s)\alpha(s)/\int_J\alpha(s)\beta(s)ds$ satisfies (\ref{normalization}),
\begin{equation*}
f_n\omega_n^\alpha=\biggl(\int_J \beta(s)\alpha(s)ds\biggr)\omega_n^\gamma
\end{equation*}
for any $n \in \bf N$. In this sense the operation of multiplication in $\mathcal T'$ has Property $\mathcal B$ (see the Introduction).
\end{remark}

\begin{example}
\label{ex3}
In the general case, let us consider the product of the delta-function $\delta_\tau^\alpha$ and a dynamic regulated function $f \in d\mathbb G$: \\ [-5mm]
\begin{equation*}
(f\delta_\tau^\alpha,\varphi)=(\delta_\tau^\alpha,f\varphi)=\int_J\varphi(\tau)(s)f(\tau)(s)\alpha(s)ds.
\end{equation*}
If $f$ has an ordinary value at $\tau$, then
the result above can be rewritten as
\begin{equation*}
f\delta_\tau^\alpha=f(\tau)\delta_\tau^\alpha.
\end{equation*}
\end{example}

\begin{remark}
The intersection $\mathcal T \cap \mathbb G$ with the topology induced by $\mathcal T$ determines a space of \textit{discontinuous test functions}. Since this is a subspace of $\mathcal T$, we obtain a factor-space of $\mathcal T'$, which is called the space of \textit{distributions with discontinuous test functions} \cite{DerKin4} (see the Introduction; also, see \cite{Kur,Kur2}).
\end{remark}

Let
$g \in s\mathbb{BV}$. Define the derivative $\dot{g} \in \mathcal T'$ by the formula
\begin{equation}
\label{dyn_deriv}
(\dot{g},\varphi)\doteq \int_I \hat{\varphi}(t)dg_c(t)+\sum_{\tau \in T(g)}\int_J\varphi(\tau)(s)(g(\tau)(s))^{\cdot}_sds.
\end{equation}
where $\varphi \in \mathcal T$, $g_c \in \mathbb{CBV}$ is a continuous part of $\ordin(g) \in \mathbb{BV}$. 
If for any $\tau \in T(g)$ we have $\sigma_\tau(g)=g(\tau+)-g(\tau-) \ne 0$, then \\ [-5mm]
\begin{equation}
\label{deriv2}
\dot{g}=\dot{g}_c+\sum_{\tau_k \in T(g)}\sigma_{\tau_k}(g)\delta_{\tau_k}^{\alpha_k}, \\ [-2mm] 
\end{equation}
where \\ [-5mm]
\begin{equation}
\label{cderiv}
(\dot{g}_c,\varphi)=\int_I \hat{\varphi}(t)dg_c(t),
\end{equation}
the shapes $\alpha_k$ are defined by
$\alpha_k(s)=(g(\tau_k)(s))^{\cdot}_s/\sigma_{\tau_k}(g)$ 
($s \in J$, $k \in \bf N$). 


\begin{theorem}
\label{dercorteo}
The functional defined by (\ref{dyn_deriv}) is linear, continuous and thus determines a distribution in $\mathcal T'$.
\end{theorem}

%

\begin{theorem}[Leibniz product rule]
\label{leibnitzteo}
For any $f,g \in s\mathbb{BV}$ we have
\begin{equation}
\label{dyn_leibnitz}
(fg)^{\cdot}=\dot{f}g+f\dot{g} \in \mathcal T'.
\end{equation}
\end{theorem}

\begin{example}
\label{dist_exdiff}
The derivative of the Heaviside function $\theta_\tau^\beta \in s\mathbb {BV}$ is given by
\begin{equation*}
\dot{\theta}_\tau^\beta=\delta_\tau^\alpha
\end{equation*}
where $\alpha=\dot{\beta}$ is a shape of the delta-function.
\end{example}


\subsection{Ordinary differential equations with distributions}

In what follows, notations $\mathcal T'_n$ and $\mathcal D_n'$
stand for the spaces of $n$-valued distributions with components in $\mathcal T'$ and $\mathcal D'$ respectively, where convergence, linear operations, operations of multiplication and differentiation are defined componentwise. We introduce analogous notations for the spaces of ordinary functions and dynamic functions $\mathbb{BV}_n$, $\mathbb{AC}_n$ and $s\mathbb{BV}_n$, respectively.

Consider in $\mathcal T_n'$
a Cauchy problem
\begin{equation}
\label{eq10}
\dot{x}=f(t,x)+g(t,x)\delta_\tau^\alpha, \quad x(t_0-)=x_0,
\end{equation}
where $t_0 \in I$, $x_0 \in \textbf{R}^n$, and the functions $f$ and $g$ with values in $\bf R^n$ and $\bf M_n$, respectively, are continuously differentiable in both variables, $\delta_\tau^\alpha \in \mathcal T_n'$ is a vector-valued delta-function,
\begin{equation*}
\alpha=(\alpha_1,\dots,\alpha_n)^{\top}, 
\end{equation*}
the value of $\delta_\tau^\alpha$ on $\varphi \in \mathcal T$ is defined by
\begin{equation}
(\delta_\tau^\alpha,\varphi)\doteq\biggl((\delta_\tau^{\alpha_1},\varphi),\dots,(\delta_\tau^{\alpha_n},\varphi)\biggr)^{\top}.
\end{equation}

We define a solution of the Cauchy problem (\ref{eq10}) to be a dynamic function $x \in s\mathbb{BV}_n$ which satisfies (\ref{eq10}) in $\mathcal T_n'$. Define an \textit{ordinary solution} $\hat{x} \doteq \ordin(x) \in \mathbb{BV}_n$.

Then there exists a solution of the Cauchy problem (\ref{eq10}) $x \in s\mathbb{BV}_n$ having ordinary values $x(t)$ for all $t \ne \tau$, and having dynamic value $x(\tau)(\cdot)$ (denote $\gamma(\cdot) \doteq x(\tau)(\cdot)$), such that
\begin{equation}
\dot{x}(t)=f\bigl(t,x(t)\bigr), 
\end{equation}
for all $t \ne \tau$, and
\begin{equation}
\label{2}
\dot{\gamma}(s)=g\bigl(\tau,\gamma(s)\bigr)\alpha(s), \quad \gamma(-1/2)=x(\tau-).
\end{equation}
(see Example \ref{ex3}), where $x(\tau+)=\gamma(1/2)$. 


Notice that in the result of the substitution of $x \in s\mathbb{BV}_n$ into (\ref{eq10}), the operations of differentiation, composition and multiplication (see example \ref{ex3}), which arise in (\ref{eq10}), are
correctly defined in $\mathcal T'_n$, in contrast to the space $\mathcal D_n'$ (see the Introduction).

\begin{remark}
Consideration of ordinary differential equations with distributions in $\mathcal D_n'$ is presented in \cite{Silv,Ramp,Ses}, where, in particular, the following Cauchy problem is considered: 
\begin{equation}
\label{bad}
\dot{x}=f(t,x)+\bigl(g(t,x)\iota\bigr)\delta_\tau, \quad x(t_0)=x_0,
\end{equation}
where $f$, $g$ are the same as above, $\delta_\tau \in \mathcal D'$ is a scalar delta-function, $\iota=(1,\dots,1)^{\top}$.

As is mentioned in \cite{Ses}, it is reasonable to expect that in the general case a solution $x$ of (\ref{bad}) is discontinuous at point $\tau$. Thus, in the general case, the notation in (\ref{bad}) is incorrect from the point of view of distribution theory, since (\ref{bad}) contains the product of a discontinuous function $g(\tau,x(\cdot))$ and a distribution $\delta_\tau \in \mathcal D'$, which is undefined in  $\mathcal D_n'$.
In \cite{Silv,Ramp,Ses} the following definition of a solution of (\ref{bad}) is proposed. The left-continuous function $x \in \mathbb{BV}_n$ is said to be a solution of (\ref{bad}), if there exists a delta-sequence $\{\omega_m\}_{m=1}^\infty$, $\omega_m \to \delta_\tau$ in $\mathcal D'$, such that 
\begin{equation}
\label{xconv}
x_m \to x
\end{equation}
($m \to \infty$) in a weak topology of $\mathbb{BV}_n$, where $x_m \in \mathbb{AC}_n$ is a solution of the approximating problem  
\begin{equation}
\dot{x}=f(t,x)+\bigl(g(t,x)\iota\bigr)\omega_m(t), \quad x(t_0)=x_0.
\end{equation}

The necessary and sufficient condition for the uniqueness of solution of (\ref{bad}), i.e., its independence on the choice of $\{\omega_m\}_{m=1}^\infty$, is a Frobenius condition
\begin{equation}
\label{frob}
\sum_{i=1}^n \frac{\partial g_{ij}}{\partial x_k}g_{km}=
\sum_{i=1}^n \frac{\partial g_{im}}{\partial x_k}g_{kj}
\end{equation}
for all $t \in I$, $x \in \bf R^n$, $k,j,m=1,\dots,n$, $g=(g_{ij})_{ij=1}^n$ \cite{Ses,Mil}. If (\ref{frob}) is satisfied,
then there exists a unique solution $x \in \mathbb{BV}_n$ of (\ref{bad}), such that 
\begin{equation}
\dot{x}=f(t,x(t)) 
\end{equation}
for all $t \ne \tau$, 
and
\begin{equation}
\label{bad2}
\dot{\eta}(s)=g(\tau,\eta(s))\iota, \quad \eta(-1/2)=x(\tau-),
\end{equation}
where $x(\tau+)=\eta(1/2)$ \cite{Ses}.

Notice that (\ref{bad2}) coincides with (\ref{2}) for $\alpha \equiv \iota$. Indeed, for (\ref{2}) the Frobenius condition (\ref{frob}) is a necessary and sufficient condition for the independence of the value $\gamma(1/2)$ on the choice of $\alpha$ satisfying (\ref{normalization}) \cite{Guyshun}.
Also, if (\ref{frob}) is satisfied, then a solution of (\ref{bad}) coincides with the ordinary solution of (\ref{eq10}).

The comparison of the results for equation (\ref{eq10}), which is considered in the space of distributions $\mathcal T_n'$, and the results for equation (\ref{bad}), which is considered in the classical space of distributions $\mathcal D_n'$, shows, that in contrast to the approach based on the space $\mathcal D_n'$ \cite{Silv,Ramp,Ses}, the use of the space of distributions $\mathcal T_n'$ allows us to make notations correct from the point of view of distribution theory, and to eliminate the restrictive Frobenius condition (\ref{frob}), which can not be satisfied in many problems of optimal control theory \cite{Mil}.
\end{remark}


\section{Multiplication of distributions by discontinuous functions in other approaches \\ [3mm]}

\subsection{Multiplication by discontinuous functions in the Schwartz space of distributions}

Various definitions of the product of a distribution and a discontinuous function in $\mathcal D'$ were proposed in \cite{Ses,Fil,Kur,Kur2,Tvr1}, and also in \cite{Sar,Sar2,Sar3}, where the operation of multiplication of a distribution by a discontinuous function is induced by the operation of multiplication of two distributions.

Furthermore, we concentrate on the multiplication of the Heaviside function $\theta_\tau$ and the delta-function. The product $\theta_\tau\delta_\tau \in \mathcal D'$ is defined in \cite{Ses,Fil,Kur,Kur2,Tvr1,Sar,Sar2,Sar3} by the formula
\begin{equation}
\label{prod100}
\theta_\tau\delta_\tau=c\delta_\tau,
\end{equation}
where $c \in \bf R$ is different in different definitions. In particular, 
in \cite{Sar,Sar2,Sar3}, where the family of the products of distributions $\{\cdot_\sigma\}_{\sigma \in \mathcal D}$ is defined (here $\mathcal D$ is a space of infinitely differentiable test functions),
for any $\sigma \in \mathcal D$ we have \\ [-6mm]
\begin{equation}
\theta_\tau \cdot_{\sigma}\delta_\tau=\sigma(0)\delta_\tau.
\end{equation}

The value $c=1/2$ is considered in \cite{Fil,Kur,Kur2}). If $c=1/2$, then we have that $(\theta_\tau \theta_\tau)^{\cdot}=\theta_\tau\delta_\tau+\delta_\tau\theta_\tau$, i.e., Leibniz product rule holds.
Also, the choice of the value $c=1/2$ is motivated by the symmetry of the delta-function $\delta_\tau \in \mathcal D'$.
Another definition of the form (\ref{prod100}) was proposed in \cite{Ses}, so the solution of the Cauchy problem (\ref{bad}) in $\mathcal D_n'$ coincides with the solution (\ref{xconv}). In \cite{Ses} the value of $c \in \bf R$ is determined by the right-hand side of the equation (\ref{bad}).

Nevertheless, any definition of the product of the form (\ref{prod100}) leads to contradictions with distribution theory, which are unavoidable in $\mathcal D'$.
Namely, let $c \in \bf R$ be given. Then we have to be able to replace delta-function in (\ref{prod100}) by the terms of arbitrary delta-sequences $\{\omega_n\}_{n=1}^\infty$, $\{\lambda_n\}_{n=1}^\infty$ in $\mathcal D'$ so that
\begin{equation}
\label{prod101}
\lim_{n \to \infty} \theta_\tau \omega_n=c\lim_{n \to \infty} \lambda_n
\end{equation}
in $\mathcal D'$. Nevertheless, for a given $c \in \bf R$ one can easily find delta-sequences $\{\omega_n\}_{n=1}^\infty$, $\{\lambda_n\}_{n=1}^\infty$ such that (\ref{prod101}) is not true. This implies that the operation of multiplication (\ref{prod100}) is not continuous in $\mathcal D'$.
Also, if $c \ne 0$, $c \ne 1$, then the operation of multiplication (\ref{prod100}) is not associative. 

In $\mathcal T'$, the product of the Heaviside function $\theta_\tau$ and the delta-function $\delta_\tau^\alpha \in \mathcal T'$ is given by the formula
\begin{equation}
\label{prod102}
\theta_\tau\delta_\tau^\alpha=\biggl(\int_0^{1/2}\alpha(s)ds\biggr)\delta_\tau^\gamma,
\end{equation}
(Example \ref{ex2}), where $\delta_\tau^\gamma \in \mathcal T'$. Let us notice the following:

1) Any definition of the form (\ref{prod100}) can be obtained from (\ref{prod102}) for the delta-functions $\delta_\tau^\alpha \in \mathcal T'$ having the shape $\alpha$ such that \\ [-5mm]
\begin{equation*}
\int_0^{1/2}\alpha(s)ds=c.
\end{equation*}

2) In contrast to (\ref{prod100}), the operation of multiplication (\ref{prod102}) is associative, satisfies the Leibniz product rule (Theorem \ref{leibnitzteo}), continuous and, furthermore, has Property $\mathcal{B}$ (see the Introduction and Remark \ref{rem1}). 

Thus, consideration of the product of the Heaviside function and the delta-function in $\mathcal T'$ allows us to avoid the aforementioned contradictions with distribution theory.

\subsection{Multiplication by discontinuous functions in the Colombeau generalized functions algebra} 

As is mentioned in the Introduction, the general multiplication of distributions is impossible in the classical space of Schwartz distributions with continuous test functions.
This defect of the Schwartz distribution theory lead to the construction of the Colombeau generalized functions algebra $\mathcal G$, containing the space of distributions $\mathcal D'$ as a subspace (here $\mathcal D'$ is a space of distributions with infinitely differentiable test functions), where such general multiplication exists \cite{Col1} (see description of the Colombeau algebra in \cite{ColM3,Obe}). The algebra of Colombeau generalized functions was applied successfully to many problems related to ordinary and partial differential equations, e.g., see \cite{Obe,ColM1,ColM2,ColM3}. 
In the classical approach the Heaviside function $\theta_\tau \in \mathcal G$ and the Dirac delta-function $\delta_\tau \in \mathcal G$ are defined as the derivatives of continuous functions in $\mathcal G$, so
$\theta_\tau$, $\delta_\tau \in \mathcal D' \subset \mathcal G$ and the product
\begin{equation*}
\theta_\tau\delta_\tau \in \mathcal G \setminus \mathcal D'
\end{equation*}
is not a distribution \cite{Col1}. However, as was observed in \cite{ColM1},\footnote{We would like to thank Professor J.F.Colombeau who kindly mentioned these results to us.} instead of consideration of the elements of the canonical embedding $\mathcal D' \subset \mathcal G$, often it is natural to work with the algebra $\mathcal G$ itself. Namely, a generalized function $f \in \mathcal G$ is said to have $\theta_\tau$ or $\delta_\tau$ as its \textit{macroscopic profile}, if 
\begin{equation*}
f \approx \theta_\tau \text{ or } f \approx \delta_\tau, 
\end{equation*}
respectively, where $\approx$ stands for the \textit{association} in $\mathcal G$. If $f$ has $\theta_\tau$ or $\delta_\tau$ as a macroscopic profile, then $f$ is called the \textit{Heaviside generalized function} or the \textit{Dirac generalized function}, respectively \cite{ColM1}. According to \cite{ColM1}, there exist different Heaviside generalized functions and Dirac generalized functions, which differ by their \textit{microscopic profiles}. The equalities of the form 
\begin{equation}
\label{colomeq}
\theta_{1,\tau}\delta_{2,\tau}=c\delta_{3,\tau}, \\ [2mm]
\end{equation}
where $\theta_{1,\tau} \in \mathcal G$ is a Heaviside generalized function,
$\delta_{2,\tau}$, $\delta_{3,\tau} \in \mathcal G$ are two Dirac generalized functions, were studied in \cite{ColM3}, where the value of constant $c \in \mathbf R$ was derived from the supplementary partial differential equation, arising from the statement of the physical problem. 
As follows from the definition of the algebra $\mathcal G$ \cite{Obe}, equality (\ref{colomeq}) can be also obtained if the Heaviside generalized function and the Dirac generalized functions in (\ref{colomeq}) are approximated by continuous functions.

The equality (\ref{colomeq}) is similar to the equality (\ref{prod_w1}) in the space of distributions $\mathcal T'$, where the constant $c \in \mathbf R$ is given explicitly (see also Remark \ref{rem1}). 
However, despite this similarity, space $\mathcal T'$ also possesses other properties. In particular, $\mathcal T'$ contains delta-functions (\ref{newdelta}) (see Remark \ref{rm_newdelta}), which allow us to unify in the subsequent paper different definitions of solution of the ordinary differential equations with distributions (see \cite{Fil} on ambiguities arising in the definition of the notion of solution of the differential equations with distributions). \\ [1mm]

\section{Proofs of statements of Sections 2 and 3 \\ [3mm]}

\begin{proof}[Proof of Lemma \ref{lem1}]
Suppose that there exists an uncountable set $A \subset I$ such that
$f(t)(\cdot)$ is not ordinary for all $t \in A$. Let us define a function $\hat{g}:I \to \bf R$ by
\begin{equation}
\label{dyn_g}
\hat{g}(t)=\sup_{s \in J} f(t)(s)-\inf_{s \in J} f(t)(s)
\end{equation}
Since $f \in d\mathbb G$ is bounded, $0 \leqslant \hat{g}(t)< \infty$ for all $t \in I$. Obviously, 
$f$ has an ordinary value at $t$ if and only if $\hat{g}(t)=0$, so $\hat{g}(t)>0$ for all $t \in A$. We define
\begin{equation*}
A_n=\biggl\{t \in A: \frac{1}{n+1}<\hat{g}(t)\leqslant \frac{1}{n}\biggr\} \quad (n=1,2,\dots),
A_0=\{t \in A: \hat{g}(t)>1\}.
\end{equation*}
Then $A=\cup_{n=0}^\infty A_n$. Since $A$ is uncountable, there exists $n_0 \geqslant 0$ such that $A_{n_0}$ is uncountable. Without loss of generality we may assume that $n_0 \ne 0$. Since $A_{n_0}$ is uncountable, there exists a closed interval $[c,d] \subset I$ such that $A_{n_0}' \doteq A_{n_0} \cap [c,d]$ is uncountable. Consequently, due to compactness of $[c,d]$, the set $A_{n_0}'$ has a limit point $\xi \in [c,d]$. Without loss of generality we may assume that there exists a right-sided neighborhood of $\xi$ which contains infinitely many points of $A_{n_0}'$. Let us show that this implies that $f(\xi+)$ does not exist.
Suppose that $p=f(\xi+) \in \bf R$ exists. Then by definition for any $\varepsilon>0$ there exists $\eta>0$ such that $\|f-p\|_{(\xi,\xi+\eta)}<\varepsilon$, i.e.,
$\sup_{t \in (a,b}\sup_{s \in J}|f(t)(s)-p|<\varepsilon$.
Now consider the value of 
$\sup_{s \in J}|f(t)(s)-p|$, 
for $t \in (\xi,\xi+\eta) \cap A_{n_0}'$,
where $(\xi,\xi+\eta) \cap A_{n_0}' \ne \varnothing$. Then the inequality 
\begin{equation}
\label{dyn_ineq1}
\sup\limits_{s \in J}|f(t)(s)-p| \geqslant 
\frac{1}{2} \biggl(\sup_{s \in J}(f(t)(s))-\inf_{s \in J}(f(t)(s)) \bigr)
\end{equation}
holds (see proof below). As follows from the definitions of $A_{n_0}'$ and $\hat{g}$, 
\begin{equation*}
\frac{1}{2}(\sup_{s \in J}(f(t)(s))-\inf_{s \in J}(f(t)(s))) \geqslant 1/n_0 
\end{equation*}
for all $t \in A_{n_0}' \cap (\xi,\xi+\eta)$. Consequently, 
$\sup_{t \in (\xi,\xi+\eta)}\sup_{s \in J}|f(t)(s)-p| \geqslant \frac{1}{2n_0}$.
Since for any $\eta>0$ $A_{n_0}' \cap (\xi,\xi+\eta) \ne \varnothing$, the last inequality holds for any $\eta>0$. This contradicts to the definition of the right-sided limit $p=f(\xi+)$.

Let us show that the inequality (\ref{dyn_ineq1}) holds. If $f(t)(s)-p \geqslant 0$ for all $s \in J$, then $\sup_{s \in J}|f(t)(s)-p|=\sup_{s \in J}(f(t)(s)-p)$ and (\ref{dyn_ineq1}) holds, where  $\int_{s \in J}(f(t)(s)-p) \geqslant 0$. In the general case we can represent $f(t)(s)=q^+(s)-q^-(s)$, where $q^+(s)$, $q^-(s) \geqslant 0$. Then
\begin{equation*}
\sup_{s \in J}|f(t)(s)-p|=\max\{\sup_{s \in J}q^+(s),\sup_{s \in J}q^-(s)\}
\geqslant 1/2(\sup_{s \in J}q^+(s)+\sup_{s \in J}q^-(s)).
\end{equation*}
Since $\sup_{s \in J}(f(t)(s)-p)=\sup_{s \in J} q^+(s)$, $-\inf_{s \in J}(f(t)(s)-p)=\sup_{s \in J} q^-(s)$, we obtain that (\ref{dyn_ineq1}) is true.
\end{proof}                                                                 

\begin{proof}[Proof of Lemma \ref{lem2}]
For a given $\tau \in I$ let us denote $p=f(\tau+)$. Then for any $\varepsilon>0$ there exits $\eta>0$ 
such that $\sup_{(\tau,\tau+\eta)}|p-f|<\varepsilon$. According to Lemma \ref{lem1}, the intersection
$(\tau,\tau+\eta) \cap \us(f) \ne \varnothing$. Then
\begin{equation*}
\sup_{t \in (\tau,\tau+\eta) \cap us(f)}|p-\hat{f}(t)|=
\sup_{t \in (\tau,\tau+\eta)\cap us(f)}\sup_{s \in J}|p-\hat{f}(t)| 
\leqslant
\sup_{(\tau,\tau+\eta)\cap us(f)}|p-f| 
\end{equation*}
where the first equality holds since $|p-\hat{f}(t)|$ does not depend on
$s \in J$. Thus, for any $t \in (\tau,\tau+\eta)$, where $\hat{f}(t)$ is defined, i.e, for all
$t \in (\tau,\tau+\eta) \cap \us(f)$, we have $|p-\hat{f}(t)|<\varepsilon$. Then the right-sided limit
$\hat{f}(\tau+)$ exists, which is equal to $f(\tau+)$. The proof for the left-sided limit is analogous.
Since $\tau \in I$ was chosen arbitrarily, the statements of the lemma holds.
\end{proof}                           

\begin{proof}[Proof of Lemma \ref{lem3}]
The proof follows from an analogous statement 
for the space of regulated functions $\mathbb G$ \cite{Hon}, the definition of a point of discontinuity of a dynamic function, and Lemma \ref{lem2}.
\end{proof}

\begin{proof}[Proof of Theorem \ref{limitteo}]                             
The linearity of the limit functional is obvious. Without loss of generality it suffices to show that $f$ is continuous at $\varphi=0$ only.
Let $\varphi_n \to 0$ in $\mathcal T$. Suppose that $(f,\varphi_n) \not\to 0$ ($n \to \infty$).
Being considering, if necessary, a subsequence of $\{(f,\varphi_n)\}_{n=1}^\infty$, we may assume that
for any $n=1,2,\dots$ the inequality $|(f,\varphi_n)|>\varepsilon_0$ holds for certain $\varepsilon_0>0$.
Due to (\ref{limit}) for any $k=1,2,\dots$ there exists $n_k$ such that $|(f_{n_k},\varphi_k)|>\frac{\varepsilon_0}{2}$.
Without loss of generality we may assume that $n_k=k$, i.e., $|(f_k,\varphi_k)|>\frac{\varepsilon_0}{2}$ for
any $k=1,2,\dots$ The last inequality contradicts to the conditions of the lemma. Thus, $(f,\varphi_n) \to 0$ as $n \to \infty$, i.e., the limit functional $f$ is continuous.
\end{proof}

\begin{proof}[Proof of Theorem \ref{ext_teo}]
According to Theorem \ref{topteo} the space of dynamic test functions $\mathcal T$ is locally-convex;
$\mathcal D$ is a subspace of $\mathcal T$. Then according to the Hahn-Banach theorem \cite{Kan},
every linear continuous functional defined on a subspace $\mathcal D$ has an extension to the whole space $\mathcal T$.
\end{proof} 

\begin{proof}[Proof of Theorem \ref{cont_teo}]
Note that $g_n \varphi \xrightarrow {} g\varphi$ in $\mathcal T$
for any $\varphi \in \mathcal T$. Consequently
\begin{multline}
\notag
|(g_nf_n,\varphi)-(gf,\varphi)|=|(f_n,g_n\varphi)-(f,g\varphi)|
\leqslant |(f_n,g_n\varphi)-(f_n,g\varphi)|+\\+|(f_n,g\varphi)-(f,g\varphi)|
\leqslant |(f_n,g_n\varphi-g\varphi)|+|(f_n,g\varphi)-(f,g\varphi)| \xrightarrow {} 0
\end{multline}
(the first and the second summands tend to zero according to Lemma \ref{lem4} 
and due to convergence $f_n \to f$ in $\mathcal T'$, respectively).
\end{proof}

\begin{proof}[Proof of Theorem \ref{dercorteo}]
The Stieltjes integral in the right-hand side of (\ref{dyn_deriv}) exists since
$\hat{\varphi} \in \mathbb G$, $\hat{g}_c \in \mathbb{CBV}$ (see \cite{Der2}), and the
convergence of the series
follows from the inequality
\begin{equation*}
\biggl|\int_J\varphi(\tau)(s)(g(\tau)(s))^{\cdot}_sds\biggr| 
\leqslant \sup_{s \in J}|\varphi(\tau)(s)|\int_J|g(\tau)(s)^{\cdot}_s|ds
\leqslant
\sup_{s \in J}|\varphi(\tau)(s)|\var_{s \in J}(g(\tau)(s))
\end{equation*}
for all $\tau \in T(g)$.
Then the following inequality holds
\begin{equation*}
\sum_{\tau \in T(g)}\biggl|\int_J\varphi(\tau)(s)(g(\tau)(s))^{\cdot}_sds\biggr| \leqslant \sup_{I}|\varphi| \\ \sum_{\tau \in T(g)}\var_{s \in J}(g(\tau)(s)),
\end{equation*}
where convergence of series in the right-hand side follows from the definition of the algebra $s\mathbb{BV}$. 
Thus, the value of $\dot{g}$ is defined for all test functions $\varphi \in \mathcal T$.
The linearity and continuity of $\dot{g}$ 
follow from definition of the convergence in $\mathcal T$, 
properties of the Stieltjes integral and properties of convergent series. 
\end{proof}


\begin{proof}[Proof of Theorem \ref{leibnitzteo}]
We have
\begin{equation}
\label{fg}
(\dot{f}g,\varphi)=(\dot{f},g\varphi)=\int_I \hat{\varphi}(t)\hat{g}(t)df_c(t)+\sum_{\tau \in T(f)}\int_J g(\tau)(s)\varphi(\tau)(s)(f(\tau)(s))^{\cdot}_s ds,
\end{equation}
\begin{equation}
\label{gf}
(f\dot{g},\varphi)=(\dot{g},f\varphi)=\int_I \hat{\varphi}(t)\hat{f}(t)dg_c(t)+\sum_{\tau \in T(g)}\int_J f(\tau)(s)\varphi(\tau)(s)(g(\tau)(s))^{\cdot}_s ds.
\end{equation}
In (\ref{fg}) and (\ref{gf}) we can perform the summation by $T(f) \cup T(g)$, since we need to add zero summands only. Consequently,
\begin{multline}
\label{fggf}
\sum_{\tau \in T(f)}\int_J g(\tau)(s)\varphi(\tau)(s)(f(\tau)(s))^{\cdot}_s ds+\sum_{\tau \in T(g)}\int_J f(\tau)(s)\varphi(\tau)(s)(g(\tau)(s))^{\cdot}_s ds=\\
=\!\sum_{\tau \in T(f) \cup T(g)}\int_J \varphi(\tau)(s)(g(\tau)(s)f(\tau)(s))^{\cdot}_sds=\!\sum_{\tau \in T(fg)}\!\int_J \varphi(\tau)(s)(g(\tau)(s)f(\tau)(s))^{\cdot}_sds,
\end{multline}
where the last equality is due to $T(fg) \subset T(f) \cup T(g)$ and $\tau \in T(f) \cup T(g) \setminus T(fg)$ if and only if $f(\tau)=0$ or $g(\tau)=0$, so we may exclude in (\ref{fggf}) the summands corresponding to $\tau \in T(f) \cup T(g) \setminus T(fg)$. We have
\begin{equation}
\label{c2}
\int_I \hat{\varphi}(t)\hat{g}(t)df_c(t)+\int_I\hat{\varphi}(t)\hat{f}(t)d_gc(t)=\int_I \hat{\varphi}(t)d\biggl(\int_a^t \hat{g}(\xi)df_c(\xi)+\int_a^t \hat{f}(\xi)dg_c(\xi) \biggr),
\end{equation}
where $I=(a,b)$.
Let us show that the following equality holds
\begin{equation}
\label{c}
(fg)_c(t)=\int_a^t \hat{g}(\xi)df_c(\xi)+\int_a^t \hat{f}(\xi)dg_c(\xi)+f(a+)g(a+)
\end{equation}
for all $t \in I$. Now $(fg)_c=(\hat{f}\hat{g})_c$, $f_c=\hat{f}_c$, $g_c=\hat{g}_c$, thus we may prove (\ref{c}) for the ordinary parts only. Hence
\begin{equation*}
\hat{f}(t)\hat{g}(t)=(\hat{f}_c(t)+\hat{f}_h(t))(\hat{g}_c(t)+\hat{g}_h(t))\!=\!
\hat{f}_c(t)\hat{g}_c(t)+\hat{f}_c(t)\hat{g}_h(t)+\hat{g}_c(t)\hat{f}_h(t)+\hat{f}_h(t)\hat{g}_h(t)
\end{equation*}
for all $t \in I$. Thus,
\begin{multline}
\notag
\int_a^t \hat{g}(\xi)d\hat{f}_c(\xi)+\int_a^t \hat{f}(\xi)d\hat{g}_c(\xi)+f(a+)g(a+)=\\=\int_a^t\hat{g}_c(\xi)d\hat{f}_c(\xi)+\int_a^t\hat{g}_h(\xi)d\hat{f}_c(\xi)+\int_a^t\hat{f}_c(\xi)d\hat{g}_c(\xi)+\int_a^t \hat{f}_h(\xi)d\hat{g}_c(\xi)=\\
=\hat{f}_c(t)\hat{g}_c(t)-\hat{f}_c(a)\hat{g}_c(a)+f(a+)g(a+)+\hat{f}_c(t)\hat{g}_h(t)-\\-\sum_{a<\tau_i<t}\hat{f}_c(\tau_i)(g(\tau_i+)-g(\tau_i-))+\hat{g}_c(t)\hat{f}_h(t)-\sum_{a<\tau_i<t}\hat{g}_c(\tau_i)(f(\tau_i+)-f(\tau_i-)),
\end{multline}
for all $t \in I$,
where $\hat{f}_c(a)\hat{f}_c(a)=f(a+)g(a+)$ and
\begin{equation*}
\hat{f}_c(t)\hat{g}_h(t)-\sum_{a<\tau_i<t}\hat{f}_c(\tau_i)(g(\tau_i+)-g(\tau_i-))
\end{equation*}
is a continuous part of $\hat{f}_c\hat{g}_h$,
\begin{equation*}
\hat{g}_c(t)\hat{f}_h(t)-\sum_{a<\tau_i<t}\hat{g}_c(\tau_i)(f(\tau_i+)-f(\tau_i-))
\end{equation*}
is a continuous part of  $\hat{g}_c\hat{f}_h$,
so equality (\ref{c}) is true.
The comparison of (\ref{fg}) -- (\ref{c}) with the definition of the derivative in $\mathcal T'$ gives us (\ref{dyn_leibnitz}).
\end{proof}

\textbf{Acknowledgments.} We would like to thank Professor J.F.~Colombeau for very important discussion and comments, Professor A.~Brudnyi for his valuable support,  Professor L.~Bates and Professor P.~Zvengrowski for their attention and help.

\bibliographystyle{amsalpha}
\bibliography{eng}

\end{document}